\documentclass[12pt,twoside]{article}
\usepackage[english]{babel}
\usepackage[latin1]{inputenc}
\usepackage{amsmath}
\usepackage{graphicx}
\usepackage{times,amssymb,amscd}
\newtheorem{thm}{Theorem}[section]

\newtheorem{lem}[thm]{Lemma}

\numberwithin{equation}{section}

\newcommand{\bA}{\mathbf{A}}

\newcommand{\bE}{\mathbf{E}}

\newcommand{\bH}{\mathbf{H}}

\newcommand{\bL}{\mathbf{L}}

\newcommand{\bR}{\mathbf{R}}
\newcommand{\bS}{\mathbf{S}}
\newcommand{\bV}{\mathbf{V}}

\newcommand{\be}{\mathbf{e}}

\newcommand{\bx}{\mathbf{x}}
\newcommand{\by}{\mathbf{y}}

\newcommand{\bT}{\mathbf{T}}

\newcommand{\bu}{\mathbf{u}}
\newcommand{\bv}{\mathbf{v}}
\newcommand{\bt}{\mathbf{t}}

\newcommand{\BV}{\boldsymbol{V}}

\newcommand{\Be}{\boldsymbol{e}}
\newcommand{\Bu}{\boldsymbol{u}}
\newcommand{\Bv}{\boldsymbol{v}}

\newcommand{\cP}{\mathcal{P}}
\newcommand{\cS}{\mathcal{S}}

\newcommand{\EUC}{\mathbf E^3}
\newcommand{\SPH}{\bS^3}
\newcommand{\HYP}{\bH^3}
\newcommand{\SXR}{\bS^2\!\times\!\bR}
\newcommand{\HXR}{\bH^2\!\times\!\bR}
\newcommand{\SLR}{\widetilde{\bS\bL_2\bR}}
\newcommand{\NIL}{\mathbf{Nil}}
\newcommand{\SOL}{\mathbf{Sol}}

\begin{document}
\pagestyle{myheadings}
\markboth{\centerline{Jen\H o Szirmai}}
{Interior angle sums of geodesic triangles $\dots$}
\title
{Interior angle sums of geodesic triangles in $\SXR$ and $\HXR$ geometries
}

\author{Jen\H o Szirmai \\
\normalsize Budapest University of Technology and \\
\normalsize Economics Institute of Mathematics, \\
\normalsize Department of Geometry \\
\normalsize Budapest, P. O. Box: 91, H-1521 \\
\normalsize szirmai@math.bme.hu
\date{\normalsize{\today}}}
\maketitle
\begin{abstract}
In the present paper we study $\SXR$ and $\HXR$ geometries, which are homogeneous Thurston 3-geometries.
We analyse the interior angle sums of geodesic triangles in both geometries and prove, 
that in $\SXR$ space it can be larger or equal than $\pi$ and in $\HXR$ space the angle sums can be less or equal than $\pi$.

In our work we will use the projective model of $\SXR$ and $\HXR$ geometries described by E. Moln\'ar in \cite{M97}.
\end{abstract}
\newtheorem{Theorem}{Theorem}[section]
\newtheorem{corollary}[Theorem]{Corollary}
\newtheorem{lemma}[Theorem]{Lemma}
\newtheorem{exmple}[Theorem]{Example}
\newtheorem{definition}[Theorem]{Definition}
\newtheorem{rmrk}[Theorem]{Remark}
\newtheorem{proposition}[Theorem]{Proposition}
\newenvironment{remark}{\begin{rmrk}\normalfont}{\end{rmrk}}
\newenvironment{example}{\begin{exmple}\normalfont}{\end{exmple}}
\newenvironment{acknowledgement}{Acknowledgement}

\section{Introduction}
\label{section1}
A geodesic triangle in Riemannian geometry and more generally in metric geometry a
figure consisting of three different points together with the pairwise-connecting geodesic curves.
The points are known as the vertices, while the geodesic curve segments are known as the sides of the triangle.

In the geometries of constant curvature $\EUC$, $\HYP$, $\SPH$ the well-known sums of the interior angles of geodesic 
triangles characterize the space. It is related to the Gauss-Bonnet theorem which states that the integral of the Gauss curvature
on a compact $2$-dimensional Riemannian manifold $M$ is equal to $2\pi\chi(M)$ where $\chi(M)$ denotes the Euler characteristic of $M$.
This theorem has a generalization to any compact even-dimensional Riemannian manifold (see e.g. \cite{Ch}, \cite{KN}).

\begin{rmrk}
In the Thurston spaces can be introduced in a natural way (see \cite{M97}) translation curves. These curves are simpler than geodesics and 
differ from them in $\NIL$, $\SLR$ and $\SOL$ geometries. In $\EUC$, $\SPH$, $\HYP$, $\SXR$ and $\HXR$ geometries the mentioned curves 
coincide with each other (\cite{B}, \cite{CsSz16}, \cite{Sz18}, \cite{Sz18-1}).
\end{rmrk}

In \cite{CsSz16} we investigated the angle sums of translation and geodesic triangles in $\SLR$ geometry
and proved that the possible sum of the interior angles in
a translation triangle must be greater or equal than $\pi$. However, in geodesic triangles this sum is
less, greater or equal to $\pi$.

In \cite{Sz16} we considered the analogous problem for geodesic triangles in $\NIL$ geometry and proved 
that the sum of the interior angles of geodesic triangles in $\NIL$ space is larger, less or equal than $\pi$.
In \cite{B} K.~Brodaczewska showed, that sum of the interior angles of translation triangles of the $\NIL$ space is larger or equal  than $\pi$.

In \cite{Sz18-1} we studied the interior angle sums of {\it translation triangles} in $\SOL$ geometry
and proved that the possible sum of the interior angles in
a translation triangle must be greater or equal than $\pi$. Further interesting properties of translation triangles and tetrahedra are described in \cite{Sz18}.

However, in $\SXR$, $\HXR$ and $\SOL$ Thur\-ston geo\-metries there are no result concerning the
angle sums of {\it geodesic triangles}. Therefore, it is interesting to study this question
in the above three geometries. 

In the present paper, we are interested in {\it geodesic triangles} in $\SXR$ and $\HXR$ spaces \cite{S,T}. 

In Section 2 we describe the projective model and the isometry group of the considered geometries,
moreover, we give an overview about its geodesic curves.
{\it In Section 3 we study the $\SXR$ and  $\HXR$ geodesic triangles and their properties.} 
\section{Projective model of $\HXR$ and $\SXR$ spaces}
E. {Moln\'ar} has shown in \cite{M97}, that the homogeneous 3-spaces
have a unified interpretation in the projective 3-sphere $\mathcal{PS}^3(\bV^4,\BV_4, \mathbf{R})$. 
In our work we shall use this projective model of $\SXR$ and $\HXR$ geometries. 
The Cartesian homogeneous coordinate simplex $E_0(\be_0)$,$E_1^{\infty}(\be_1)$,$E_2^{\infty}(\be_2)$,
$E_3^{\infty}(\be_3)$, $(\{\be_i\}\subset \bV^4$ \ $\text{with the unit point}$ $E(\be = \be_0 + \be_1 + \be_2 + \be_3 ))$ 
which is distinguished by an origin $E_0$ and by the ideal points of coordinate axes, respectively. 
Moreover, $\by=c\bx$ with $0<c\in \mathbf{R}$ (or $c\in\mathbf{R}\setminus\{0\})$
defines a point $(\bx)=(\by)$ of the projective 3-sphere $\cP \cS^3$ (or that of the projective space $\cP^3$ where opposite rays
$(\bx)$ and $(-\bx)$ are identified). 
The dual system $\{(\Be^i)\}\subset \BV_4$ describes the simplex planes, especially the plane at infinity 
$(\Be^0)=E_1^{\infty}E_2^{\infty}E_3^{\infty}$, and generally, $\Bv=\Bu\frac{1}{c}$ defines a plane $(\Bu)=(\Bv)$ of $\cP \cS^3$
(or that of $\cP^3$). Thus $0=\bx\Bu=\by\Bv$ defines the incidence of point $(\bx)=(\by)$ and plane
$(\Bu)=(\Bv)$, as $(\bx) \text{I} (\Bu)$ also denotes it. Thus {$\SXR$} can be visualized in the affine 3-space $\bA^3$
(so in $\bE^3$) as well.
\subsection{Geodesic curves in $\SXR$ space}
In this section we recall the important notions and results from the papers \cite{M97}, \cite{PSSz10}, \cite{Sz13-1}, \cite{Sz11-1}, \cite{Sz11-2}. 

The well-known infinitezimal arc-length square at any point of $\SXR$ as follows
\begin{equation}
   \begin{gathered}
     (ds)^2=\frac{(dx)^2+(dy)^2+(dz)^2}{x^2+y^2+z^2}.
       \end{gathered} \tag{2.1}
     \end{equation}
We shall apply the usual geographical coordiantes $(\phi, \theta), ~ (-\pi < \phi \le \pi, ~ -\frac{\pi}{2}\le \theta \le \frac{\pi}{2})$ 
of the sphere with the fibre coordinate $t \in \bR$. We describe points in the above coordinate system in our model by the following equations: 
\begin{equation}
x^0=1, \ \ x^1=e^t \cos{\phi} \cos{\theta},  \ \ x^2=e^t \sin{\phi} \cos{\theta},  \ \ x^3=e^t \sin{\theta} \tag{2.2}.
\end{equation}
Then we have $x=\frac{x^1}{x^0}=x^1$, $y=\frac{x^2}{x^0}=x^2$, $z=\frac{x^3}{x^0}=x^3$, i.e. the usual Cartesian coordinates.
We obtain by \cite{M97} that in this parametrization the infinitezimal arc-length square 
at any point of $\SXR$ is the following
\begin{equation}
   \begin{gathered}
      (ds)^2=(dt)^2+(d\phi)^2 \cos^2 \theta +(d\theta)^2.
       \end{gathered} \tag{2.3}
     \end{equation}
The geodesic curves of $\SXR$ are generally defined as having locally minimal arc length between their any two (near enough) points. 
The equation systems of the parametrized geodesic curves $\gamma(t(\tau),\phi(\tau),\theta(\tau))$ in our model can be determined by the 
general theory of Riemann geometry (see \cite{KN}, \cite{Sz11-2}).

Then by (2.2) we get with $c=\sin{v}$, $\omega=\cos{v}$ the equation systems of a geodesic curve, visualized in Fig.~3 in our Euclidean model:
\begin{equation}
  \begin{gathered}
   x(\tau)=e^{\tau \sin{v}} \cos{(\tau \cos{v})}, \\ 
   y(\tau)=e^{\tau \sin{v}} \sin{(\tau \cos{v})} \cos{u}, \\
   z(\tau)=e^{\tau \sin{v}} \sin{(\tau \cos{v})} \sin{u},\\
   -\pi < u \le \pi,\ \ -\frac{\pi}{2}\le v \le \frac{\pi}{2}. \tag{2.4}
  \end{gathered}
\end{equation}
\begin{definition}
The distance $d(P_1,P_2)$ between the points $P_1$ and $P_2$ is defined by the arc length of the shortest geodesic curve 
from $P_1$ to $P_2$.
\end{definition}
\subsection{Geodesic curves of $\HXR$ geometry}
In this section we recall the important notions and results from the papers \cite{M97},  \cite{PSSz11}, \cite{Sz12-1}.

The points of $\HXR$ space, forming an open cone solid in the projective space $\mathcal{P}^3$, are the following:
\begin{equation}
\HXR:=\big\{ X(\bx=x^i \be_i)\in \mathcal{P}^3: -(x^1)^2+(x^2)^2+(x^3)^2<0<x^0,~x^1 \big\}. \notag
\end{equation}
In this context E. Moln\'ar \cite{M97} has derived the infinitezimal arc-length square at any point of $\HXR$ as follows
\begin{equation}
   \begin{gathered}
     (ds)^2=\frac{1}{(-x^2+y^2+z^2)^2}\cdot [(x)^2+(y)^2+(z)^2](dx)^2+ \\ + 2dxdy(-2xy)+2dxdz (-2xz)+ [(x)^2+(y)^2-(z)^2] (dy)^2+ \\ 
     +2dydz(2yz)+ [(x)^2-(y)^2+(z)^2](dz)^2.
       \end{gathered} \tag{2.5}
     \end{equation}
This becomes simpler in the following special (cylindrical) coordiantes $(t, r, \alpha)$, $(r \ge 0, ~ -\pi < \alpha \le \pi)$ 
with the fibre coordinate $t \in \bR$. We describe points in our model by the following equations: 
\begin{equation}
x^0=1, \ \ x^1=e^t \cosh{r},  \ \ x^2=e^t \sinh{r} \cos{\alpha},  \ \ x^3=e^t \sinh{r} \sin{\alpha}  \tag{2.6}.
\end{equation}
Then we have $x=\frac{x^1}{x^0}=x^1$, $y=\frac{x^2}{x^0}=x^2$, $z=\frac{x^3}{x^0}=x^3$, i.e. the usual Cartesian coordinates.
We obtain by \cite{M97} that in this parametrization the infinitezimal arc-length square by (2.1)
at any point of $\HXR$ is the following
\begin{equation}
   \begin{gathered}
      (ds)^2=(dt)^2+(dr)^2 +\sinh^2{r}(d\alpha)^2.
       \end{gathered} \tag{2.7}
     \end{equation}
The geodesic curves of $\HXR$ are generally defined as having locally minimal arc length between their any two (near enough) points. 
The equation systems of the parametrized geodesic curves $\gamma(t(\tau),r(\tau),\alpha(\tau))$ in our model can be determined by the 
general theory of Riemann geometry:

By (2.5) the second order differential equation system of the $\HXR$ geodesic curve is the following \cite{Sz12-1}:
\begin{equation}
\ddot{\alpha}+2\coth(r) ~ \dot{r} \dot{\alpha}=0, ~ \ddot{r}-\sinh(r) \cosh(r) \dot{\alpha}^2=0, ~ \ddot{t}=0, \tag{2.8}
\end{equation}
from which we get first a line as "geodesic hyperbola" on our model of $\bH^2$ times a component on $\bR$ each running with constant velocity $c$ and $\omega$, respectively:
\begin{equation}
t=c\cdot \tau, \ \ \alpha=0, \ \ r=\omega \cdot \tau, \ \ c^2+\omega^2=1 \tag{2.9}.
\end{equation}
We can assume, that the starting point of a geodesic curve is $(1,1,0,0)$, because we can transform a curve into an 
arbitrary starting point, moreover, unit velocity with "geographic" coordinates $(u,v)$ can be assumed:
\begin{equation}
\begin{gathered}
        r(0)=\alpha(0)=t(0)=0; \ \ \dot{t}(0)= \sin{v}, \ \dot{r}(0)=\cos{v} \cos{u}, \dot{\alpha}(0)=\cos{v} \sin{u}; \\
        - \pi < u \leq \pi, ~ -\frac{\pi}{2}\le v \le \frac{\pi}{2}. \notag
\end{gathered}
\end{equation}
Then by (2.6) we get with $c=\sin{v}$, $\omega=\cos{v}$ the equation systems of a geodesic curve, visualized in Fig.~8 in our Euclidean model \cite{Sz12-1}:
\begin{equation}
  \begin{gathered}
   x(\tau)=e^{\tau \sin{v}} \cosh{(\tau \cos{v})}, \\ 
   y(\tau)=e^{\tau \sin{v}} \sinh{(\tau \cos{v})} \cos{u}, \\
   z(\tau)=e^{\tau \sin{v}} \sinh{(\tau \cos{v})} \sin{u},\\
   -\pi < u \le \pi,\ \ -\frac{\pi}{2}\le v \le \frac{\pi}{2}. \tag{2.10}
  \end{gathered}
\end{equation}
\begin{definition}
The distance $d(P_1,P_2)$ between the points $P_1$ and $P_2$ is defined by the arc length of the geodesic curve 
from $P_1$ to $P_2$.
\end{definition}
\begin{rmrk}
$\SXR$ and $\HXR$ are affine metric spaces (affine-projective spaces -- in the sense of the unified formulation of \cite{M97}). Therefore their linear, affine, unimodular,
etc. transformations are defined as those of the embedding affine space.
\end{rmrk}
\section{Geodesic triangles}
We consider $3$ points $A_1$, $A_2$, $A_3$ in the projective model of $X$ space (see Section 2) $(X\in\{\SXR, \HXR \}$.
The {\it geodesic segments} $a_k$ connecting the points $A_i$ and $A_j$
$(i<j,~i,j,k \in \{1,2,3\}, k \ne i,j$) are called sides of the {\it geodesic triangle} with vertices $A_1$, $A_2$, $A_3$ (see Fig.~1,~2).

In Riemannian geometries the infinitesimal arc-lenght square (see (2.1) and (2.5)) is used to define the angle $\theta$ between two geodesic curves.
If their tangent vectors in their common point are $\bu$ and $\bv$ and $g_{ij}$ are the components of the metric tensor then
\begin{equation}
\cos(\theta)=\frac{u^i g_{ij} v^j}{\sqrt{u^i g_{ij} u^j~ v^i g_{ij} v^j}} \tag{3.1}
\end{equation}
It is clear by the above definition of the angles and by the infinitesimal arc-lenght squares that
the angles are the same as the Euclidean ones at the starting point of the geodesics.

Considering a geodesic triangle $A_1A_2A_3$ we can assume by the homogeneity of the considered geometries that one of its vertex 
coincide with the point $A_1=(1,1,0,0)$ and the other two vertices are $A_2=(1,x^2,y^2,z^2)$ and $A_3=(1,x^3,y^3,z^3)$. 

We will consider the {\it interior angles} of geodesic triangles that are denoted at the vertex $A_i$ by $\omega_i$ $(i\in\{1,2,3\})$.
We note here that the angle of two intersecting geodesic curves depends on the orientation of their tangent vectors. 
\subsection{Interior angle sums in $\SXR$ geometry}

In order to determine the interior angles of a geodesic triangle $A_1A_2A_3$ 
and its interior angle sum $\sum_{i=1}^3(\omega_i)$,
we define {\it isometric transformations} $\bT^{\SXR}_{A_i}$, $(i\in \{2,3\}$, as elements of the isometry group of $\SXR$ geometry, that
maps the $A_i$ onto $A_1$).
Let the isometrie $\bT^{\SXR}_{A_2}$ be given by the composition of some special types of $\SXR$ isometries, which transforms 
a fixed $A_2=(1,x_2,y_2,z_2)$ point of $\SXR$ into $(1,1,0,0)$ (up to a positive determinant factor):

${\mathcal{T}}=({\bf Id. }, T)$ is a fibre translation, 
\begin{equation}
\begin{gathered}
A_2=(1,x_2,y_2,z_2)\rightarrow A_2^{\mathcal{T}}=(1,x_2',y_2',z_2')= \\=
A_2^{\mathcal{T}}=\Big(1,\frac{x_2}{\sqrt{x_2^2+y_2^2+z_2^2}},\frac{y_2}{\sqrt{x_2^2+y_2^2+z_2^2}},
\frac{z_2}{\sqrt{x_2^2+y_2^2+z_2^2}}\Big).
\end{gathered} \tag{3.2}
\end{equation}
($A_2^{{\mathcal{T}}}$ has $0$ fibre coordinate).
${\mathcal{R}_x}=({\bf R_x},0)$ is a special rotation about $x$ axis with $0$ fibre translation, which moves the point $(1,x_2',y_2',z_2')$
into the $[x,y]$ plane.
\begin{equation}
\begin{gathered}
A_2^{{\mathcal{T}}}=(1,x_2',y_2',z_2') \rightarrow A_2^{{\mathcal{T}}{\mathcal{R}_x}}=(1,x_2'',y_2'',0)=\\=
A_2^{{\mathcal{T}}{\mathcal{R}_x}}=(1,x_2',\sqrt{y_2'^2+z_2'^2},0).
\end{gathered} \tag{3.3}
\end{equation}  
Similarly, ${\mathcal{R}_z}=({\bf R_z},0)$ is a special rotation about $z$ axis with $0$ fibre translation, 
which moves the point $(1,x_2'',y_2'',0)$ into the $(1,1,0,0)$ point.
\begin{equation}
\begin{gathered}
A_2^{\mathcal{T}{\mathcal{R}_x}}=(1,x_2'',y_2'',0) \rightarrow A_2^{{\mathcal{T}}{\mathcal{R}_x}{\mathcal{R}_z}}=(1,1,0,0).
\end{gathered} \tag{3.4}
\end{equation}  
Finally we apply the inverse transformation $\mathcal{R}^{-1}_x$ of rotation $\mathcal{R}_x$ because of that the geodesic curve $g(A_1,A_2)$ between the points 
$A_1$ and $A_2$ and its image $g(A_1^2,A_1)$
under the transformation ${\mathcal{T}}{\mathcal{R}_x}{\mathcal{R}_z}{\mathcal{R}^{-1}_x}$ lie in the same plane in Euclidean sense.
The matrix of the above transformation $\bT^{\SXR}_{A_2}={\mathcal{T}}{\mathcal{R}_x}{\mathcal{R}_z}{\mathcal{R}^{-1}_x}$ is the following:
\begin{figure}[ht]
\centering
\includegraphics[width=18cm]{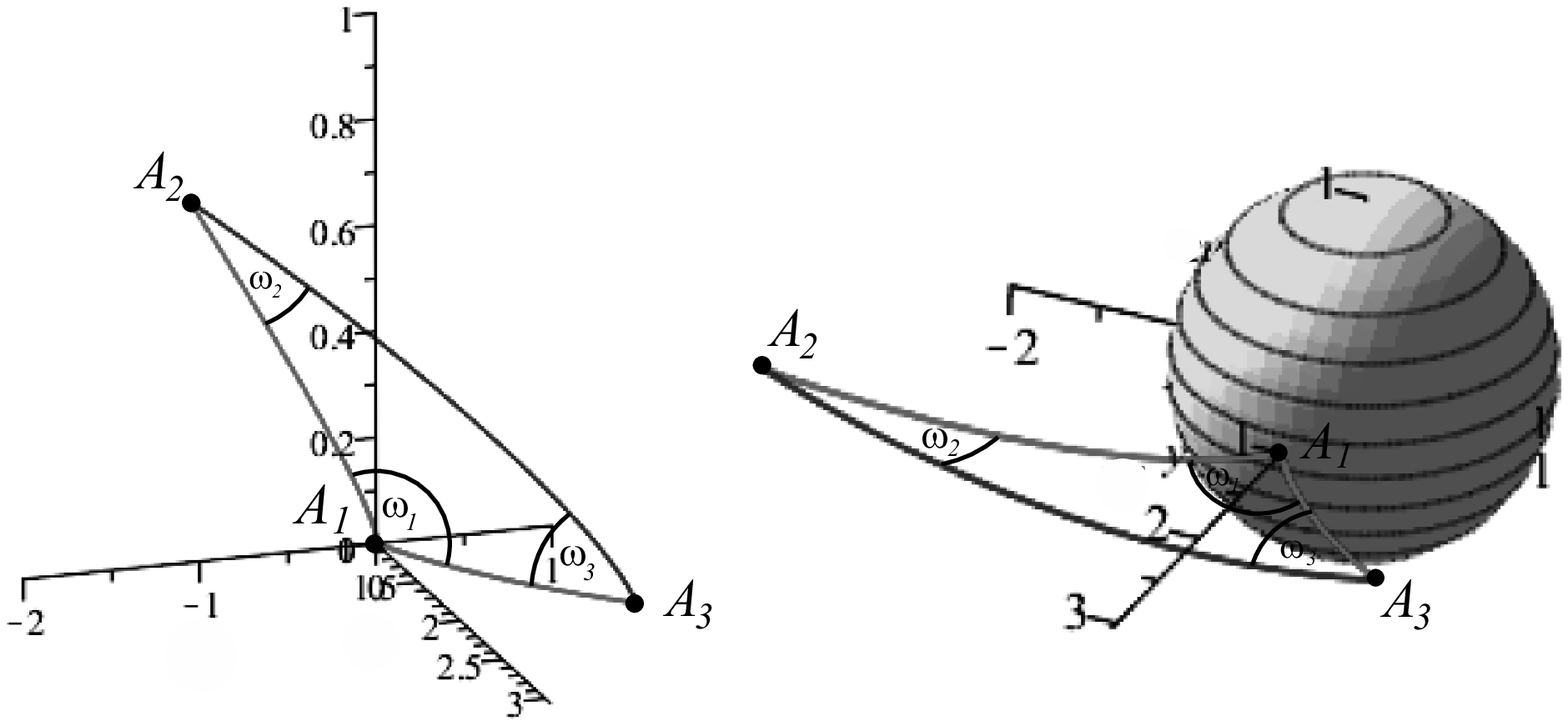}
\caption{Geodesic triangle with vertices $A_1=(1,1,0,0)$, $A_2=(1,3,-2,1)$, $A_3=(1,2,1,0)$ in $\SXR$ geometry.}
\label{}
\end{figure}

$\bT^{\SXR}_{A_2}=$ 
{\footnotesize 
\begin{equation}
\begin{gathered}
=\begin{pmatrix}
1 & 0 & 0 & 0 \\
0 & \frac{x_2}{(x_2)^2+(y_2)^2+(z_2)^2} & \frac{-y_2}{(x_2)^2+(y_2)^2+(z_2)^2} & \frac{-z_2}{(x_2)^2+(y_2)^2+(z_2)^2} \\
0 & \frac{y_2}{(x_2)^2+(y_2)^2+(z_2)^2} & \frac{(y_2)^2 x_2+(z_2)^2\sqrt{(x_2)^2+(y_2)^2+(z_2)^2}}{((x_2)^2+(y_2)^2+(z_2)^2)((y_2)^2+(z_2)^2)} & \frac{-y_2 z_2(-x_2 + \sqrt{(x_2)^2+(y_2)^2+(z_2)^2})}{((x_2)^2+(y_2)^2+(z_2)^2)((y_2)^2+(z_2)^2)}  \\
0 & \frac{z_2}{(x_2)^2+(y_2)^2+(z_2)^2} & \frac{(y_2) z_2(-x_2+\sqrt{(x_2)^2+(y_2)^2+(z_2)^2})}{((x_2)^2+(y_2)^2+(z_2)^2)((y_2)^2+(z_2)^2)} & \frac{(z_2)^2 x_2+ (y_2)^2 \sqrt{(x_2)^2+(y_2)^2+(z_2)^2})}{((x_2)^2+(y_2)^2+(z_2)^2)((y_2)^2+(z_2)^2)} 
\end{pmatrix} ,  \tag{3.5}
\end{gathered}
\end{equation}}
and the images $\bT^{\SXR}_{A_2}(A_i)$ of the vertices $A_i$ $(i \in \{1,2,3\})$ are the following (see also Fig.~2):
\begin{equation}
\begin{gathered}
\bT^{\SXR}_{A_2}(A_1)=A_1^2=\\ =\Big(1,\frac{x_2}{(x_2)^2+(y_2)^2+(z_2)^2},\frac{-y_2}{(x_2)^2+(y_2)^2+(z_2)^2},\frac{-z_2}{(x_2)^2+(y_2)^2+(z_2)^2}\Big),\\
\bT^{\SXR}_{A_2}(A_2)=A_2^2=(1,1,0,0), \\ \bT^{\SXR}_{A_2}(A_3)=
A_3^2=\Big(1,\frac{x_2x_3+y_2y_3}{(x_2)^2+(y_2)^2+(z_2)^2}, \\
\frac{y_3(z_2)^2\sqrt{(x_2)^2+(y_2)^2+(z_2)^2}+x_2(y_2)^2y_3-x_3(y_2)^3-x_3y_2(z_2)^2}{((y_2)^2+(z_2)^2)((x_2)^2+(y_2)^2+(z_2)^2)},\\
-\frac{z_2(y_3 y_2(\sqrt{(x_2)^2+(y_2)^2+(z_2)^2}-x_2)+x_3(y_2)^2+x_3(z_2)^2)}{((y_2)^2+(z_2)^2)((x_2)^2+(y_2)^2+(z_2)^2)}\Big). \tag{3.6}
\end{gathered}
\end{equation}
\begin{rmrk}
More informations about the isometry group of $\SXR$ and about its discrete subgroups can be found in \cite{Sz11-1} and \cite{Sz11-2}.  
\end{rmrk}
Similarly to the above computation we get that the images $\bT^{\SXR}_{A_3}(A_i)$ of the vertices $A_i$ $(i \in \{1,2,3\})$ are the following (see also Fig.~2):
\begin{equation}
\begin{gathered}
\bT^{\SXR}_{A_3}(A_1)=A_1^3=\Big(1,\frac{x_3}{(x_3)^2+(y_3)^2},\frac{-y_3}{(x_3)^2+(y_3)^2},0\Big),\\
\bT^{\SXR}_{A_3}(A_3)=A_3^3=A_1=(1,1,0,0), \\ \bT^{\SXR}_{A_3}(A_2)=
A_2^3=\Big(1,\frac{x_2x_3+y_2y_3}{(x_3)^2+(y_3)^2}, 
\frac{x_3y_2-x_2y_3}{(x_3)^2+(y_3)^2}, 
\frac{z_2}{\sqrt{(x_3)^2+(y_3)^2}}\Big). \tag{3.7}
\end{gathered}
\end{equation}
\begin{figure}[ht]
\centering
\includegraphics[width=13cm]{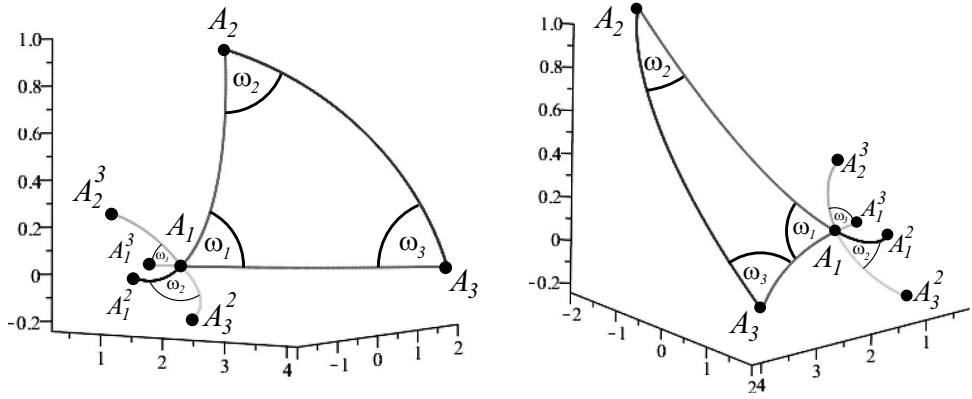}
\caption{Geodesic triangle with vertices $A_1=(1,1,0,0)$, $A_2=(1,3,-2,1)$, $A_3=(1,2,1,0)$ in $\SXR$ geometry, and transformed images of its geodesic side segments.}
\label{}
\end{figure}
Our aim is to determine angle sum $\sum_{i=1}^3(\omega_i)$ of the interior angles of geodesic triangles $A_1A_2A_3$ (see Fig.~1,~2).
We have seen that $\omega_1$ and the angle of geodesic curves with common point at the vertex $A_1$ is the same as the 
Euclidean one therefore can be determined by usual Euclidean sense.

The $\bT^{\SXR}_{A_i}$ $(i=2,3)$ are isometries in $\SXR$ geometry thus
$\omega_i$ is equal to the angle $(g(A_i^i, A_1^i),g(A_i^i, A_j^i))\angle$ $(i,j=2,3$, $i \ne j)$ (see Fig.~2)
where $g(A_i^i, A_1^i)$, $g(A_i^i, A_j^i)$ are oriented geodesic curves $(A_1=A_2^2=A_3^3)$ and
$\omega_1$ is equal to the angle $(g(A_1, A_2),g(A_1, A_3)) \angle$ 
where $g(A_1, A_2)$, $g(A_1, A_3)$ are also oriented geodesic curves.

We denote the oriented unit tangent vectors of the oriented geodesic curves $g(A_1, A_i^j)$ with $\mathbf{t}_i^j$ where
$(i,j)\in\{(1,3),(1,2),(2,3),(3,2),(3,0),(2,0)\}$ and $A_3^0=A_3$, $A_2^0=A_2$.
The Euclidean coordinates of $\mathbf{t}_i^j$ (see Section 2.1) are :
\begin{equation}
\mathbf{t}_i^j=(\sin(v_i^j), \cos(v_i^j) \cos(u_i^j), \cos(v_i^j) \sin(u_i^j)). \tag{3.8}
\end{equation}
In order to obtain the angle of two geodesic curves $g{(A_1,A_i^j})$ and $g{(A_1,A_k^l)}$ ($(i,j)\ne(k,l)$; $(i,j),(k,l)\in\{(1,3),(1,2),(2,3),(3,2),(3,0),(2,0)\})$ 
intersected at the vertex $A_1$ we need to determine their tangent vectors $\bt_s^r$ 
$((s,r) \in \{(1,3),(1,2),$ $(2,3),(3,2),(3,0),(2,0)\})$
(see (3.8)) at their starting point $A_1$.
From (3.8) follows that a tangent vector at the origin is given by the parameters $u$ and $v$ of the corresponding geodesic curve (see (2.10)) that 
can be determined from the homogeneous coordinates of the endpoint of the geodesic curve as the following Lemma shows:
\begin{lem}
Let $(1,x,y,z)$ $(x,y,z \in \bR, x^2+y^2+z^2 \ne 0)$ be the homogeneous coordinates of the point $P \in \SXR$. The paramerters of the 
corresponding geodesic curve $g{(A_1,P)}$ are the following:
\begin{enumerate}
\item $y,z \in \bR \setminus \{0\}$ and $x^2+y^2+z^2 \ne 1$;
\begin{equation}
\begin{gathered}
v=\mathrm{arctan}\Big(\frac{\log \sqrt{x^2+y^2+z^2}}{\mathrm{arccos}\frac{x}{\sqrt{x^2+y^2+z^2}}}\Big),~u=\mathrm{arctan}\Big(\frac{z}{y}\Big),\\
\tau=\frac{\log \sqrt{x^2+y^2+z^2}}{\sin v}, ~ \text{where} ~ -\pi < u \le \pi, ~ -\pi/2\le v \le \pi/2, ~ \tau \in \bR^+.
\end{gathered} \tag{3.9}
\end{equation} 
\item $y=0$, $z\ne 0$ and $x^2+z^2 \ne 1$;
\begin{equation}
\begin{gathered}
u=\frac{\pi}{2}, ~v=\mathrm{arctan}\Big(\frac{\log \sqrt{x^2+z^2}}{\mathrm{arccos}\frac{x}{\sqrt{x^2+z^2}}}\Big),\\
\tau=\frac{\log \sqrt{x^2+z^2}}{\sin v}, ~ \text{where}  ~ -\pi/2\le v \le \pi/2, ~ \tau \in \bR^+.
\end{gathered} \tag{3.10}
\end{equation}
\item $y=0$, $z\ne 0$ and $x^2+z^2 = 1$;
\begin{equation}
\begin{gathered}
u=\frac{\pi}{2}, ~v=0, ~ \tau=\arccos(x), ~ \tau \in \bR^+.
\end{gathered} \tag{3.11}
\end{equation}
\item $y, z =0$;
\begin{equation}
u=0,~v=\frac{\pi}{2},~\tau=\log \sqrt{x^2+y^2+z^2}, ~ \tau \in \bR^+. \tag{3.12} 
\end{equation}
\item $x=0,~y=0$ and $z \ne 1$;
\begin{equation}
\begin{gathered}
u=\frac{\pi}{2},~ v=\mathrm{arctan}\frac{2  \log|z|}{\pi}, ~\tau=\frac{\log|z|}{\sin v}, \\ -\pi/2\le v \le \pi/2, ~ \tau \in \bR^+. \tag{3.13}
\end{gathered}
\end{equation}
\end{enumerate} ~ ~ $\square$
\end{lem}
\begin{figure}[ht]
\centering
\includegraphics[width=12cm]{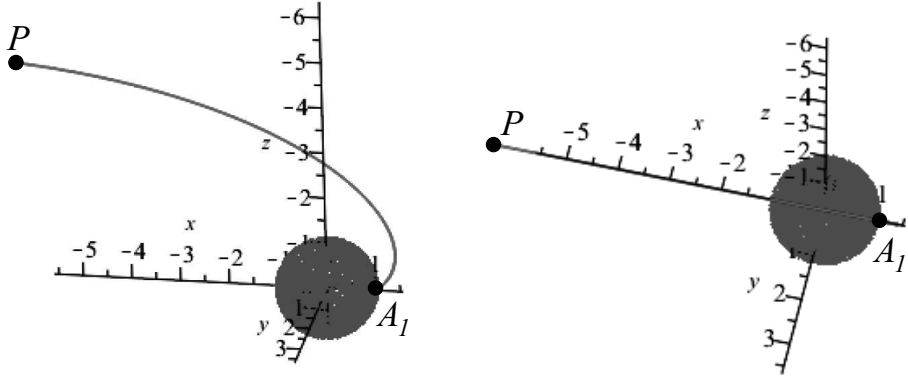}
\caption{Geodesic curve $g(A_1,P)$ ($A_1=(1,1,0,0)$ and $P \in \SXR$) with ``base plane", the plane of a geodesic curve contain the origin $E_0=(1,0,0,0)$ of the model.}
\label{}
\end{figure}
We obtain directly from the (2.4) equations of the geodesic curves the following
\begin{lem}
Let $P$ be an arbitrary point and $g(A_1,P)$ ($A_1=(1,1,0,0)$) is a geodesic curve in the considered model of $\SXR$ geometry. 
The points of the geodesic curve $g(A_1,P)$ and the centre of the model $E_0$ lie in a plane in Euclidean sense (see Fig.~3). ~ ~ $\square$
\end{lem}
\begin{thm}
If the Euclidean plane of the vertices of a $\SXR$ geodesic triangle $A_1A_2A_3$ contains the centre of model $E_0$ then its 
interior angle sum is equal to $\pi$.
\end{thm}
\begin{figure}[ht]
\centering
\includegraphics[width=13cm]{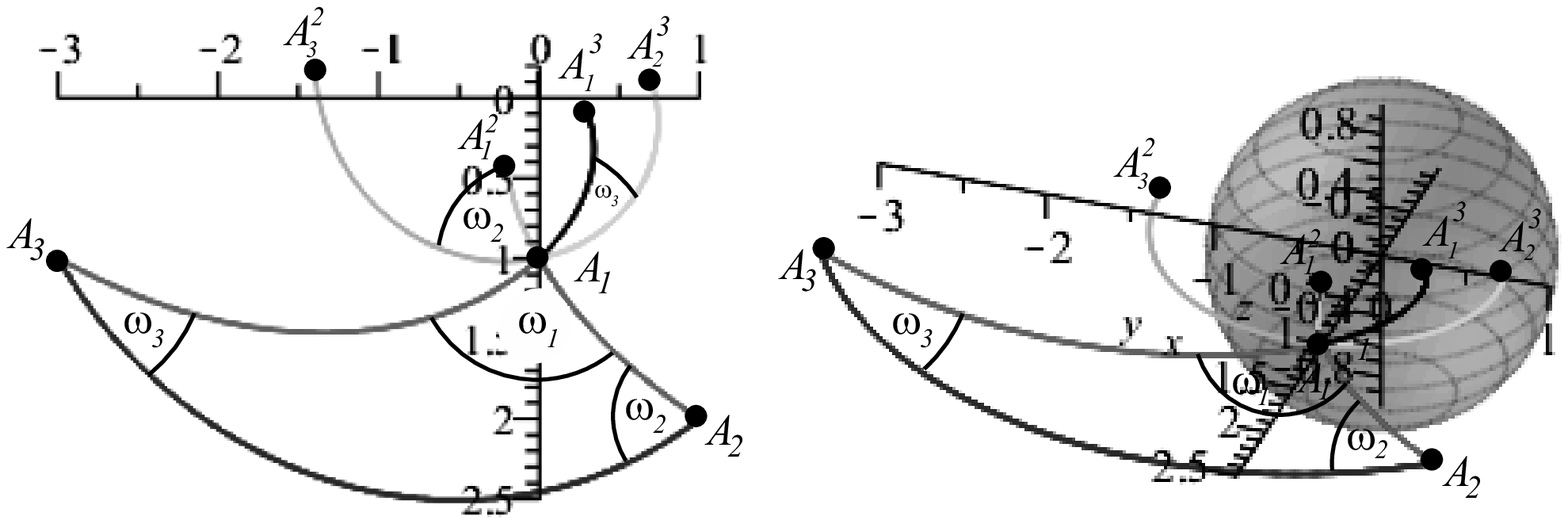}
\caption{Geodesic triangle with vertices $A_1=(1,1,0,0)$, $A_2=(1,1,-3,0)$, $A_3=(1,2,1,0)$ in $\SXR$ geometry, and transformed images of its geodesic side segments.
The geodesic curve segments $g{(A_1,A_2)}$, $g{(A_2,A_3)}$, $g({A_1,A_3})$ lie
on the coordinate plane $[x,y]$ and the interior angle sum of this geodesic triangle is $\sum_{i=1}^3(\omega_i)=\pi$.}
\label{}
\end{figure}
\textbf{Proof:} We can assume without loss of generality that the vertices $A_1,A_2,A_3$ of such a geodesic 
triangle lie in the $[x,y]$ plane of the model. Using the Lemma 3.2 we get, that the geodesic segments $A_1A_2$, $A_1A_3$ and $A_2A_3$ are containing by the [x,y] plane, too.

The $\SXR$ transformations $\bT_{A_2}^{\SXR}$ and $\bT_{A_3}^{\SXR}$ are isometries 
in $\SXR$ geometry thus $\omega_2$ is equal to the angle $(g(A_2^2, A_1^2), g(A_2^2, A_3^2)) \angle$ (see Fig.~2,~4) 
of the oriented geodesic segments $g({A_2^2, A_1^2})$, $g({A_2^2,A_3^2})$ and $\omega_3$ is equal to the angle 
$(g(A_3^3, A_1^3),g(A_3^3 A_2^3)) \angle$ 
of the oriented geodesic segments $g({A_3^3, A_1^3})$ and $g({A_3^3, A_2^3})$ $(A_1=A_2^2=A_3^3$).  

Substituting the coordinates of the points $A_i^j$ (see (3.5), (3.6) and (3.7)) $((i,j) \in \{(1,3),(1,2),$ $(2,3),(3,2),(3,0),(2,0)\})$ to the appropriate equations (3.8-12) of Lemma 3.1, 
it is easy to see that 
\begin{equation}
\begin{gathered}
v_2^0=-v_1^2,~u_2^0-u_1^2=\pm \pi \Rightarrow \bt_2^0=-\bt_1^2,\\ 
v_3^0=-v_1^3,~u_3^0-u_1^3=\pm \pi \Rightarrow \bt_3^0=-\bt_1^3,\\
v_3^2=-v_2^3,~u_3^2-u_2^3=\pm \pi \Rightarrow \bt_3^2=-\bt_2^3.
\end{gathered} \tag{3.13}
\end{equation}
The endpoints $T_i^j$ of the position vectors $\bt_i^j=\overrightarrow{A_1T_i^j}$  
lie on the unit sphere centred at the origin. The measure of angle $\omega_i$ $(i\in \{1,2,3\})$ of the vectors $\bt_i^j$ and $\bt_r^s$ is equal to the spherical 
distance of the corresponding points $T_i^j$ and $T_r^s$ on the unit sphere (see Fig.~4). Moreover, a direct consequence of equations (3.13) that each point pair 
($T_2$, $T_1^2$), $(T_3$,$T_1^3$), ($T_2^3$,$T_3^2$) 
contains antipodal points related to the unit sphere with centre $A_1$.

Due to the antipodality $\omega_1=T_2A_1T_3 \angle =T_1^2A_1T_1^3 \angle$, therefore their corresponding spherical 
distances are equal, as well (see Fig.~4). 
Now, the sum of the interior angles $\sum_{i=1}^3(\omega_i)$ can be considered as three consecutive spherical arcs $(T_3^2 T_1^2)$, $(T_1^2 T_1^3)$, 
$T_1^3 T_2^3)$. 
Since the points $T_2$, $T_1^2$, $T_3$, $T_1^3$, $T_2^3$, $T_3^2$ lie in the $[x,y]$ plane (see Lemma 3.2) the sum of these arc lengths is equal to the half 
of the circumference of the main circle on the unit sphere {i.e.} $\pi$. $\square$
\medbreak 
\begin{figure}[ht]
\centering
\includegraphics[width=10cm]{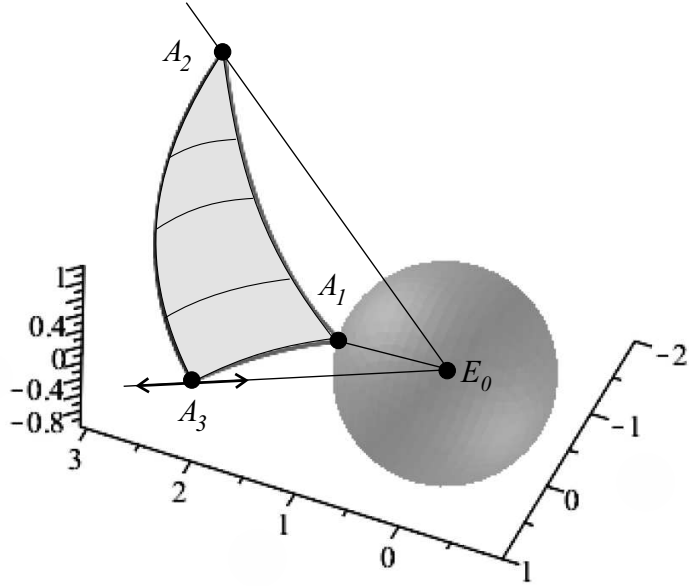}
\caption{Geodesic triangle with vertices $A_1=(1,1,0,0)$, $A_2=(1,3,-2,1)$, $A_3=(1,2,1,0)$ and the correspondig trihedron with base sphere of $\SXR$ geometry.}
\label{}
\end{figure}
\begin{figure}[ht]
\centering
\includegraphics[width=7cm]{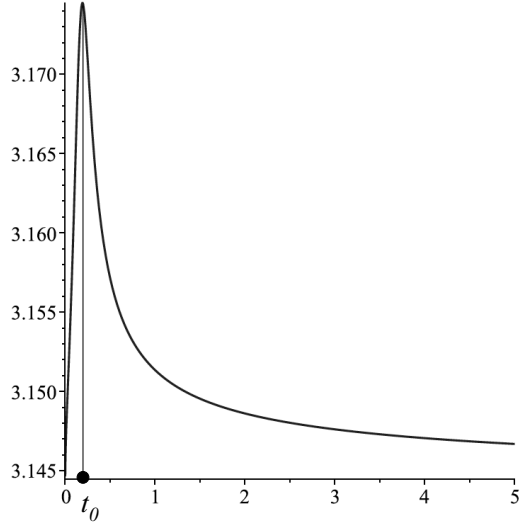}
\caption{$S(\Delta(t))$ function related to parameters $x_2=3, y_2=-2, z_2=1$ $x_3=2\cdot t,y_3=1 \cdot t,z_3=0$.}
\label{}
\end{figure}
We can determine the interior angle sum of arbitrary geodesic triangle.
In the following table we summarize some numerical data of interior angles of given geodesic triangles:
\medbreak
\centerline{\vbox{
\halign{\strut\vrule\quad \hfil $#$ \hfil\quad\vrule
&\quad \hfil $#$ \hfil\quad\vrule &\quad \hfil $#$ \hfil\quad\vrule &\quad \hfil $#$ \hfil\quad\vrule &\quad \hfil $#$ \hfil\quad\vrule
\cr
\noalign{\hrule}
\multispan5{\strut\vrule\hfill\bf Table 1: ~ $A_1=(1,0,0,0)$, $A_2=(1,3,-2,1)$  \hfill\vrule}%
\cr
\noalign{\hrule}
\noalign{\vskip2pt}
\noalign{\hrule}
A_3/\omega_i & \omega_1 & \omega_2 & \omega_3  & \sum_{i=1}^3(\omega_i)  \cr
\noalign{\hrule}
(1,2/\sqrt{5},1/\sqrt{5},0) & 1.97206 & 0.26028 &  0.92635 & 3.15869 \cr
\noalign{\hrule}
(1,2,1,0) & 0.94654  & 0.68775 & 1.51707 &  3.15135 \cr
\noalign{\hrule}
(1,4,2,0) & 0.73193 & 1.29546 & 1.12123 & 3.14862 \cr
\noalign{\hrule}
(1,12,6,0) & 0.61470 & 1.99926 & 0.53246 &  3.14643 \cr
\noalign{\hrule}
(1,2000,1000,0) & 0.50628 & 2.52677 & 0.11050  & 3.14355 \cr
\noalign{\hrule}
}}}
\medbreak
\medbreak

By the above experiences and computations we obtain the following 
\begin{thm}
If the Euclidean plane of the vertices of a $\SXR$ geodesic triangle $A_1A_2A_3$ does not contain the centre of model $E_0$ then its 
interior angle sum is greater than $\pi$.
\end{thm}
\textbf{Proof:} We can assume without loss of generality that the vertices $A_1,A_2$ of such a geodesic 
triangle lie in the $[x,y]$ plane of the model. Using the Lemma 3.2 we get, that the geodesic segment $A_iA_j$, ($(i,j)\in\{(1,2),(1,3),2,3)\}$) is contained by the $A_iA_jE_0$ plane,
therefore the sides of triangle $A_1A_2A_3$ lie in the boundary of trihedron given by the points $E_0$, $A_1$, $A_2$, $A_3$ (see Fig.~2 and 5). It is clear, that all types of geodesic triangles
can be described by such a triangle. 
Therefore, it is sufficient investigate the interior angle sums of geodesic triangles where we fix two of the vertices, e.g. $A_1$ and $A_2$ and move the third vertex 
$A_3$ on the half straight line $E_0A_3$ with starting point $E_0 \ne A_3(t)$. 
\begin{rmrk}
It is well known, that if the vertices $A_1,A_2,A_3$ lie in a sphere of radius $R\in \bR^+$ centred at $E_0$ then the interior angle sum of spherical triangle 
$A_1A_2A_3$ is greater than $\pi$.
\end{rmrk}
Let $\Delta^{\SXR}(t)$ $(t\in\bR^+)$ denote the above geodesic triangle with {\it interior angles} at 
the vertex $A_i$ by $\omega_i(t)$ $(i\in\{1,2,3\})$.  

The interior angle sum function $S(\Delta(t))=\sum_{i=1}^3(\omega_i(t))$ can be determined related to the parameters
$x_2,y_2,z_2,x_3,y_3 \in \bR$ by the formulas (2.4), (3.6), (3.7) and by the Lemma 3.1.
Analyzing the above complicated continuous functions of single real variable $t$ we get that its maximum is achieved at a point $t_0 \in (0, \infty)$ depending on given parameters. 
Moreover, $S(\Delta^{\SXR}(t))$ is stricly increasing on the interval $(0,t_0)$, stricly decreasing on the interval $(t_0,\infty)$ and
$$
\lim_{t \rightarrow 0}S(\Delta^{\SXR}(t))=\pi, ~ ~ ~ ~ \lim_{t \rightarrow \infty} S(\Delta^{\SXR}(t))=\pi. \hspace{1cm} 
$$ 
In Fig.~6 we described the $S(\Delta^{\SXR}(t))$ function related to geodesic triangle $\Delta^{\SXR}(t)$ $(t \in (0,5))$ with vertices $A_1=(1,1,0,0)$, $A_2=(1,3,-2,1)$, $A_3=(1,2\cdot t,1 \cdot t,0)$.
Its maximum is achieved at $t_0 \approx 0.19316$ where $S(\Delta^{\SXR}(t_0))\approx 3.17450$. ~ ~ ~ $\square$

Finally we get the following
\begin{thm}
The sum of the interior angles of a geodesic triangle of $\SXR$ space is greater or equal to $\pi$. ~ ~ $\square$
\end{thm}
\subsection{Interior angle sums in $\HXR$ geometry}
Similarly to the $\SXR$ space we investigate the interior angles of a geodesic triangle $A_1A_2A_3$ 
and its interior angle sum $\sum_{i=1}^3(\omega_i)$ in the $\HXR$ space.
Therefore we define {\it isometric transformations} $\bT^{\HXR}_{A_i}$, $(i\in \{2,3\})$, as elements of the isometry group of $\HXR$ geometry, that
maps the $A_i$ onto the vertex $A_1$.
Let the isometrie $\bT^{\HXR}_{A_2}$ be given by the composition of some special types of $\HXR$ isometries, which transforms 
a fixed $A_2=(1,x_2,y_2,z_2)$ point of $\HXR$ into $A_1=(1,1,0,0)$ (up to a positive determinant factor). 
The methods, the considered transformations and the determinations of their matrices are similar to the
$\SXR$ case and are therefore not detailed here.
\begin{figure}[ht]
\centering
\includegraphics[width=13cm]{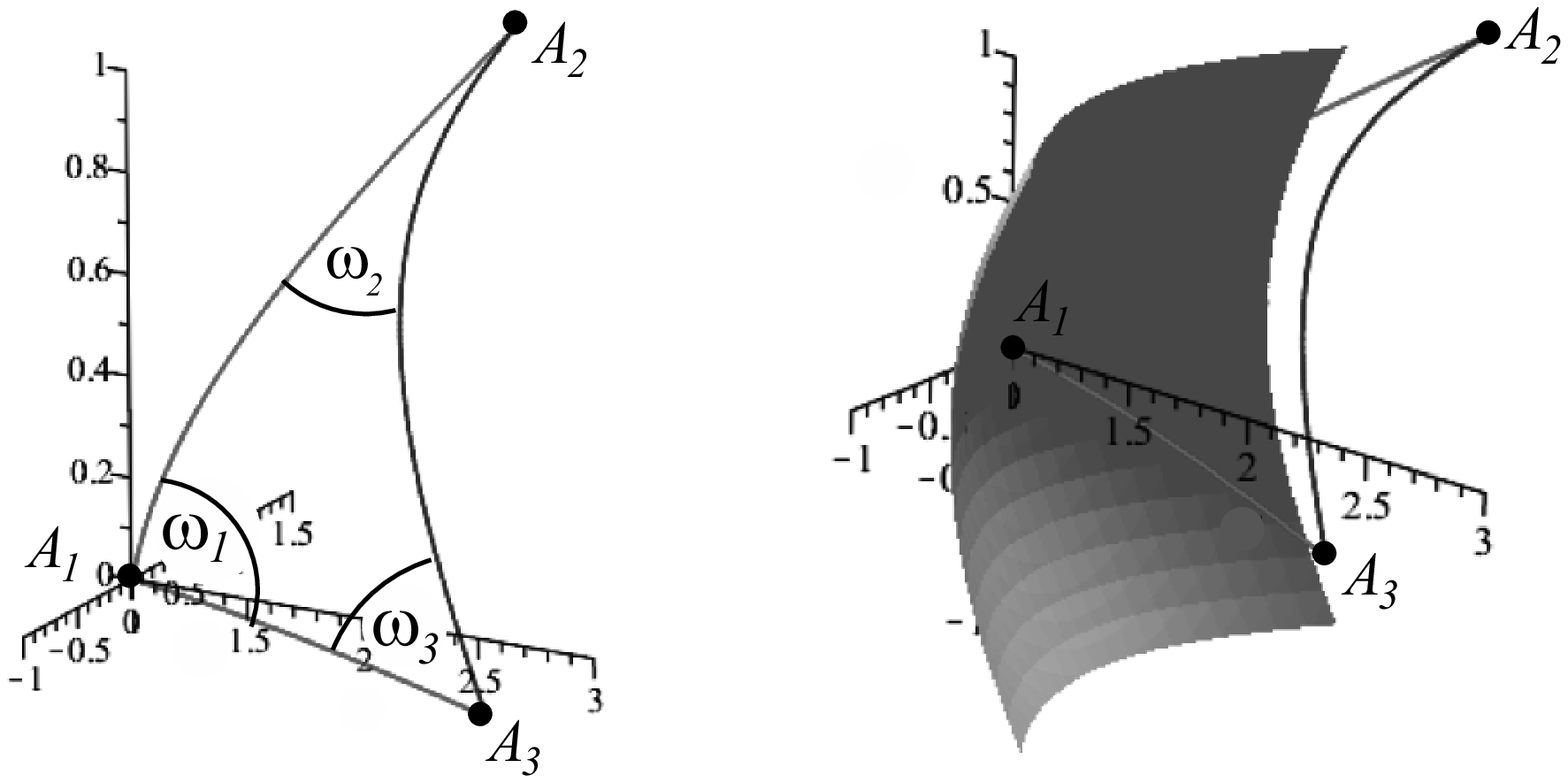}
\caption{Geodesic triangle with vertices $A_1=(1,1,0,0)$, $A_2=(1,2,3/2,1)$, $A_3=(1,3,-1,0)$ in $\HXR$ geometry.}
\label{}
\end{figure}
The images $\bT^{\HXR}_{A_2}(A_i)$ of the vertices $A_i$ $(i \in \{1,2,3\})$ are the following (see also Fig.~7,~9):
\begin{equation}
\begin{gathered}
\bT^{\HXR}_{A_2}(A_1)=A_1^2=\\ =\Big(1,\frac{x_2}{(x_2)^2-(y_2)^2-(z_2)^2},\frac{-y_2}{(x_2)^2-(y_2)^2-(z_2)^2},\frac{-z_2}{(x_2)^2-(y_2)^2-(z_2)^2}\Big),\\
\bT^{\HXR}_{A_2}(A_2)=A_2^2=(1,1,0,0), \\ \bT^{\HXR}_{A_2}(A_3)=
A_3^2=\Big(1,\frac{x_2x_3-y_2y_3}{(x_2)^2-(y_2)^2-(z_2)^2}, \\
\frac{y_3(z_2)^2\sqrt{(x_2)^2-(y_2)^2-(z_2)^2}+x_2(y_2)^2y_3-x_3(y_2)^3-x_3y_2(z_2)^2}{((y_2)^2+(z_2)^2)((x_2)^2-(y_2)^2-(z_2)^2)},\\
-\frac{z_2(y_3 y_2(\sqrt{(x_2)^2-(y_2)^2-(z_2)^2}-x_2)+x_3(y_2)^2+x_3(z_2)^2)}{((y_2)^2+(z_2)^2)((x_2)^2-(y_2)^2-(z_2)^2)}\Big). \tag{3.14}
\end{gathered}
\end{equation}
\begin{rmrk}
More informations about the isometry group of $\HXR$ and about its discrete subgroups can be found in \cite{Sz12-1}.  
\end{rmrk}
Similarly to the above computation we get that the images $\bT^{\HXR}_{A_3}(A_i)$ of the vertices $A_i$ $(i \in \{1,2,3\})$ are the following (see also Fig.~7,~9):
\begin{equation}
\begin{gathered}
\bT^{\HXR}_{A_3}(A_1)=A_1^3=\Big(1,\frac{x_3}{(x_3)^2-(y_3)^2},\frac{-y_3}{(x_3)^2-(y_3)^2},0\Big),\\
\bT^{\HXR}_{A_3}(A_3)=A_3^3=A_1=(1,1,0,0), \\ \bT^{\HXR}_{A_3}(A_2)=
A_2^3=\Big(1,\frac{x_2x_3-y_2y_3}{(x_3)^2-(y_3)^2}, 
\frac{x_3y_2-x_2y_3}{(x_3)^2-(y_3)^2}, 
\frac{z_2}{\sqrt{(x_3)^2-(y_3)^2}}\Big). \tag{3.15}
\end{gathered}
\end{equation}
The method is the same as that used for $\SXR$ case to determine angle sum $\sum_{i=1}^3(\omega_i)$ of 
the interior angles of geodesic triangles $A_1A_2A_3$ (see Fig.~7,~9).
We have seen that $\omega_1$ and the angle of geodesic curves with 
common point at the vertex $A_1$ is the same as the 
Euclidean one therefore can be determined by usual Euclidean sense.

$\omega_i$ is equal to the angle $(g(A_i^i, A_1^i),g(A_i^i, A_j^i))\angle$ $(i,j=2,3$, $i \ne j)$ (see Fig.~7,~9)
where $g(A_i^i, A_1^i)$, $g(A_i^i, A_j^i)$ are oriented geodesic curves $(A_1=A_2^2=A_3^3)$ and
$\omega_1$ is equal to the angle $(g(A_1, A_2),g(A_1, A_3)) \angle$ 
where $g(A_1, A_2)$, $g(A_1, A_3)$ are also oriented geodesic curves.
\begin{figure}[ht]
\centering
\includegraphics[width=12cm]{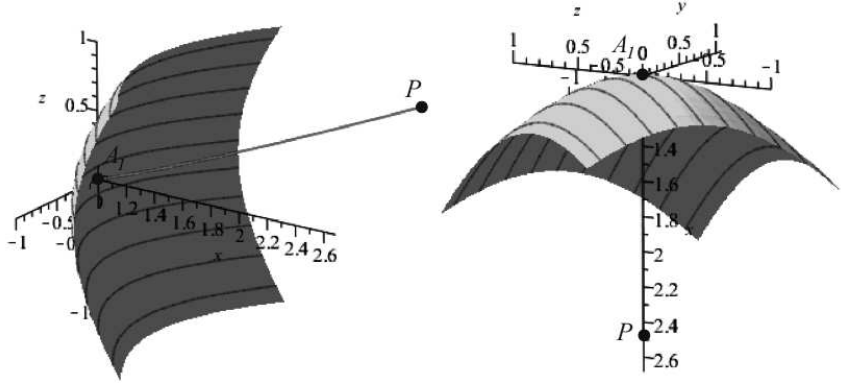}
\caption{Geodesic curve $g(A_1,P)$ ($A_1=(1,1,0,0)$ and $P \in \HXR$) with ``base plane" (the "upper" sheet of the two-sheeted hyperboloid), the plane of a geodesic curve contain the origin $E_0=(1,0,0,0)$ of the model.}
\label{}
\end{figure}
We denote the oriented unit tangent vectors of the oriented geodesic curves $g(A_1, A_i^j)$ with $\mathbf{t}_i^j$ where
$(i,j)\in\{(1,3),(1,2),$ $(2,3),(3,2),(3,0),(2,0)\}$ and $A_3^0=A_3$, $A_2^0=A_2$.
The Euclidean coordinates of $\mathbf{t}_i^j$ coincide with the coordinates in (3.8) 
(see Section 2.2).
In order to obtain the angle of two geodesic curves $g{(A_1,A_i^j})$ and $g{(A_1,A_k^l)}$ ($(i,j)\ne(k,l)$; $(i,j),(k,l)\in\{(1,3),(1,2),$ $(2,3),(3,2),(3,0),(2,0)\})$ 
intersected at the vertex $A_1$ we need to determine their tangent vectors $\bt_s^r$ 
$((s,r) \in \{(1,3),(1,2),$ $(2,3),(3,2),(3,0),(2,0)\})$
(see (2.10) and (3.8)) at their starting point $A_1$. From (3.8) follows that a tangent vector at the origin is given by the parameters $u$ and $v$ of the 
corresponding geodesic curve (see (2.10)) that 
can be determined from the homogeneous coordinates of the endpoint of the geodesic curve as the following Lemma shows:
\begin{lem}
Let $(1,x,y,z)$ $(x,y,z \in \bR, x^2-y^2-z^2 \ge 0,~ x \ge 0)$ be the homogeneous coordinates of the point $P \in \HXR$. The paramerters of the 
corresponding geodesic curve $g{(A_1,P)}$ are the following:
\begin{enumerate}
\item $y,z \in \bR \setminus \{0\}$ and $x^2-y^2-z^2 \ne 1$;
\begin{equation}
\begin{gathered}
v=\mathrm{arctan}\Big(\frac{\log \sqrt{x^2-y^2-z^2}}{\mathrm{arccosh}\frac{x}{\sqrt{x^2-y^2-z^2}}}\Big),~u=\mathrm{arctan}\Big(\frac{z}{y}\Big),\\
\tau=\frac{\log \sqrt{x^2-y^2-z^2}}{\sin v}, ~ \text{where} ~ -\pi < u \le \pi, ~ -\pi/2\le v \le \pi/2, ~ \tau \in \bR^+.
\end{gathered} \tag{3.16}
\end{equation} 
\item $y=0$, $z\ne 0$ and $x^2-z^2 \ne 1$;
\begin{equation}
\begin{gathered}
u=\frac{\pi}{2}, ~v=\mathrm{arctan}\Big(\frac{\log \sqrt{x^2-z^2}}{\mathrm{arccosh}\frac{x}{\sqrt{x^2-z^2}}}\Big),\\
\tau=\frac{\log \sqrt{x^2-z^2}}{\sin v}, ~ \text{where}  ~ -\pi/2\le v \le \pi/2, ~ \tau \in \bR^+.
\end{gathered} \tag{3.17}
\end{equation}
\item $y=0$, $z\ne 0$ and $x^2-z^2 = 1$;
\begin{equation}
\begin{gathered}
u=\frac{\pi}{2}, ~v=0, ~ \tau=\mathrm{arccosh}(x), ~ \tau \in \bR^+.
\end{gathered} \tag{3.18}
\end{equation}
\item $y, z =0$;
\begin{equation}
u=0,~v=\frac{\pi}{2},~\tau=\log (x), ~ \tau \in \bR^+. \tag{3.19} 
\end{equation}
\end{enumerate} ~ ~ $\square$
\end{lem}
We obtain directly from the (2.10) equations of the geodesic curves the following
\begin{lem}
Let $P$ be an arbitrary point and $g(A_1,P)$ ($A_1=(1,1,0,0)$) is a geodesic curve in the considered model of $\HXR$ geometry. 
The points of the geodesic curve $g(A_1,P)$ and the centre of the model $E_0$ lie in a plane in Euclidean sense (see Fig.~8).~ ~ $\square$
\end{lem}
The proof of the next theorem  essentially is the same as the proof of Theorem 3.3.
\begin{thm}
If the Euclidean plane of the vertices of a $\HXR$ geodesic triangle $A_1A_2A_3$ contains the centre of model $E_0$ then its 
interior angle sum is equal to $\pi$ (see Fig.~9). ~ ~ $\square$
\end{thm}
\begin{figure}[ht]
\centering
\includegraphics[width=13cm]{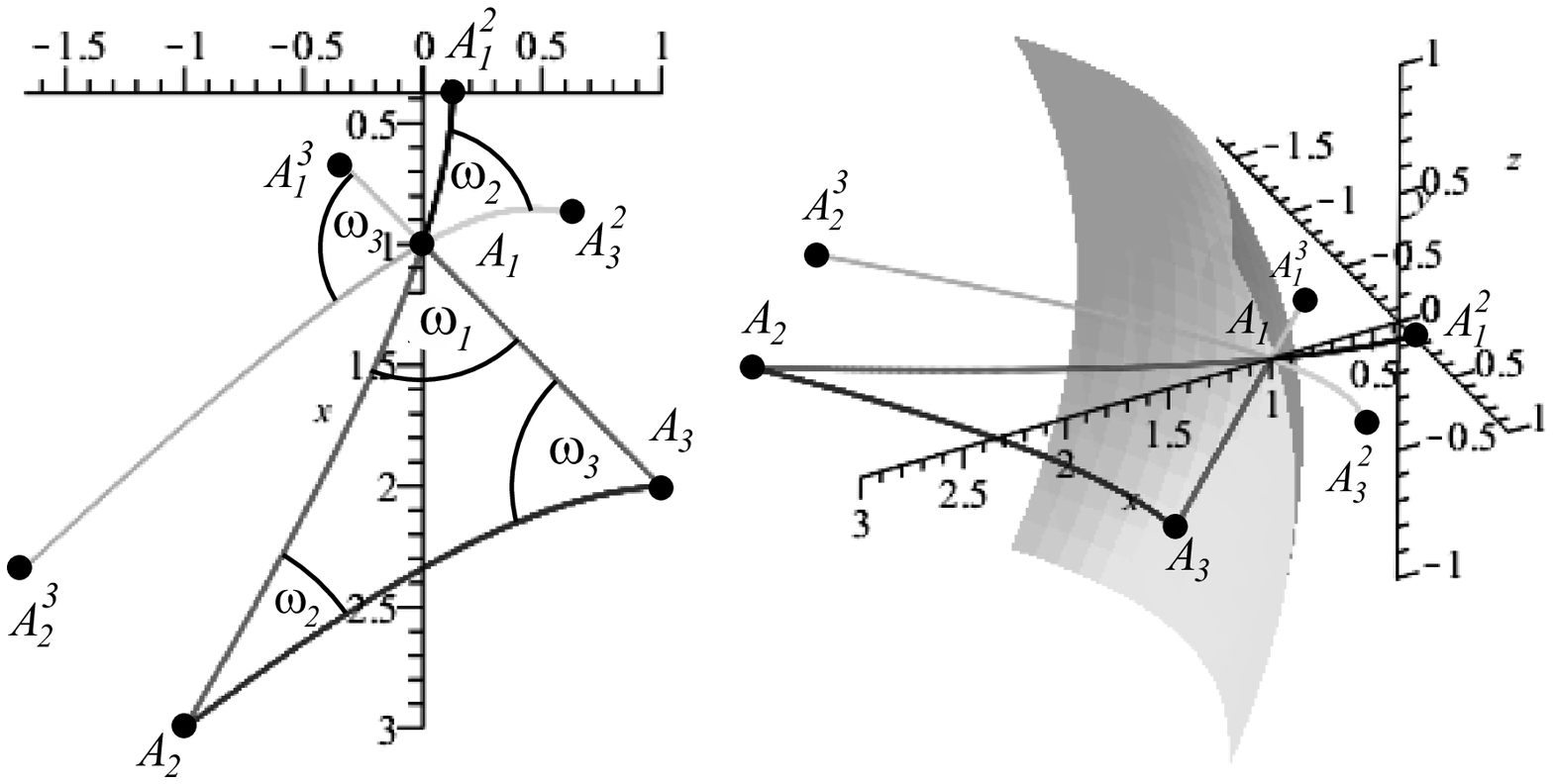}
\caption{Geodesic triangle with vertices $A_1=(1,1,0,0)$, $A_2=(1,2,3/2,1)$, $A_3=(1,3,-1,0)$ in $\HXR$ geometry, and transformed images of its geodesic side segments.
The geodesic curve segments $g{(A_1,A_2)}$, $g{(A_2,A_3)}$, $g({A_3,A_1})$ lie
on the coordinate plane $[x,y]$ and the interior angle sum of this geodesic triangle is $\sum_{i=1}^3(\omega_i)=\pi$.}
\label{}
\end{figure}
We can determine the interior angle sum of arbitrary $\HXR$ geodesic triangle.
In the following table we summarize some numerical data of interior angles of given geodesic triangles:

\medbreak
\centerline{\vbox{
\halign{\strut\vrule\quad \hfil $#$ \hfil\quad\vrule
&\quad \hfil $#$ \hfil\quad\vrule &\quad \hfil $#$ \hfil\quad\vrule &\quad \hfil $#$ \hfil\quad\vrule &\quad \hfil $#$ \hfil\quad\vrule
\cr
\noalign{\hrule}
\multispan5{\strut\vrule\hfill\bf Table 2: ~ $A_1=(1,0,0,0)$, $A_2=(1,2,3/2,1)$ \hfill\vrule}%
\cr
\noalign{\hrule}
\noalign{\vskip2pt}
\noalign{\hrule}
A_3 & \omega_1 & \omega_2 & \omega_3  & \sum_{i=1}^3(\omega_i)  \cr
\noalign{\hrule}
(1,3/\sqrt{8},-1/\sqrt{8},0) & 2.54659 & 0.06953 & 0.41780 & 3.03392 \cr
\noalign{\hrule}
(1,3,-1,0) & 1.93230  & 0.49280 & 0.69816 & 3.12325 \cr
\noalign{\hrule}
(1,6,-2,0) & 1.83102 & 0.71611 & 0.58348 & 3.13061 \cr
\noalign{\hrule}
(1,9,-3,0) & 1.80083  & 0.81224 & 0.51964 &  3.13270 \cr
\noalign{\hrule}
(1,3000,-1000,0) & 1.70394 & 1.25735 & 0.17793 & 3.13922 \cr
\noalign{\hrule}
}}}
By the above experiences and computations we obtain the following 
\begin{thm}
If the Euclidean plane of the vertices of a $\HXR$ geodesic triangle $A_1A_2A_3$ does not contain the centre of model $E_0$ then its 
interior angle sum is less than $\pi$.
\end{thm}
{\bf Proof:} The proof is similar to the $\SXR$ case. 

We can assume without loss of generality that the vertices $A_1,A_2$ of such a geodesic 
triangle lie in the $[x,y]$ plane of the model. Using the Lemma 3.7 we get, that the geodesic segment $A_iA_j$, ($(i,j)\in\{(1,2),(1,3),2,3)\}$) is contained by the $A_iA_jE_0$ plane,
therefore the sides of triangle $A_1A_2A_3$ lie in the boundary of trihedron given by the points $E_0$, $A_1$, $A_2$, $A_3$. It is clear, that all types of geodesic triangles
can be described by such a triangle. 
Therefore, it is sufficient investigate the interior angle sums of geodesic triangles where we fix two of the vertices, e.g. $A_1$ and $A_2$ and move the third vertex 
$A_3$ on the half straight line $E_0A_3$ with starting point $E_0 \ne A_3(t)$. 
\begin{rmrk}
It is well known, that if the vertices $A_1,A_2,A_3$ lie in a "upper" sheet of the two-sheeted hyperboloid (in the hyperboloid model of the hyperbolic plane geometry where 
the straight lines of hyperbolic 2-space are modeled by geodesics on the hyperboloid) centred at $E_0$ then the interior angle sum of hyperbolic triangle 
$A_1A_2A_3$ is less than $\pi$.
\end{rmrk}
Let $\Delta(t)$ $(t\in\bR^+)$ denote the above geodesic triangle with {\it interior angles} at 
the vertex $A_i$ by $\omega_i(t)$ $(i\in\{1,2,3\})$.  

The interior angle sum function $S(\Delta^{\HXR}(t))=\sum_{i=1}^3(\omega_i(t))$ can be determined related to the parameters
$x_2,y_2,z_2,x_3,y_3 \in \bR$ by the formulas (2.10), (3.14), (3.15) and by the Lemma 3.6.
Analyzing the above complicated continuous functions of single real variable $t$ we get that its maximum is achieved at a point $t_0 \in (0, \infty)$ depending on given parameters. 
Moreover, $S(\Delta^{\HXR}(t))$ is stricly increasing on the interval $(0,t_0)$, stricly decreasing on the interval $(t_0,\infty)$ and
$$
\lim_{t \rightarrow 0}S(\Delta^{\HXR}(t))=\pi, ~ ~ ~ ~ \lim_{t \rightarrow \infty} S(\Delta^{\HXR}(t))=\pi. \hspace{1cm} 
$$ 
\begin{figure}[ht]
\centering
\includegraphics[width=7cm]{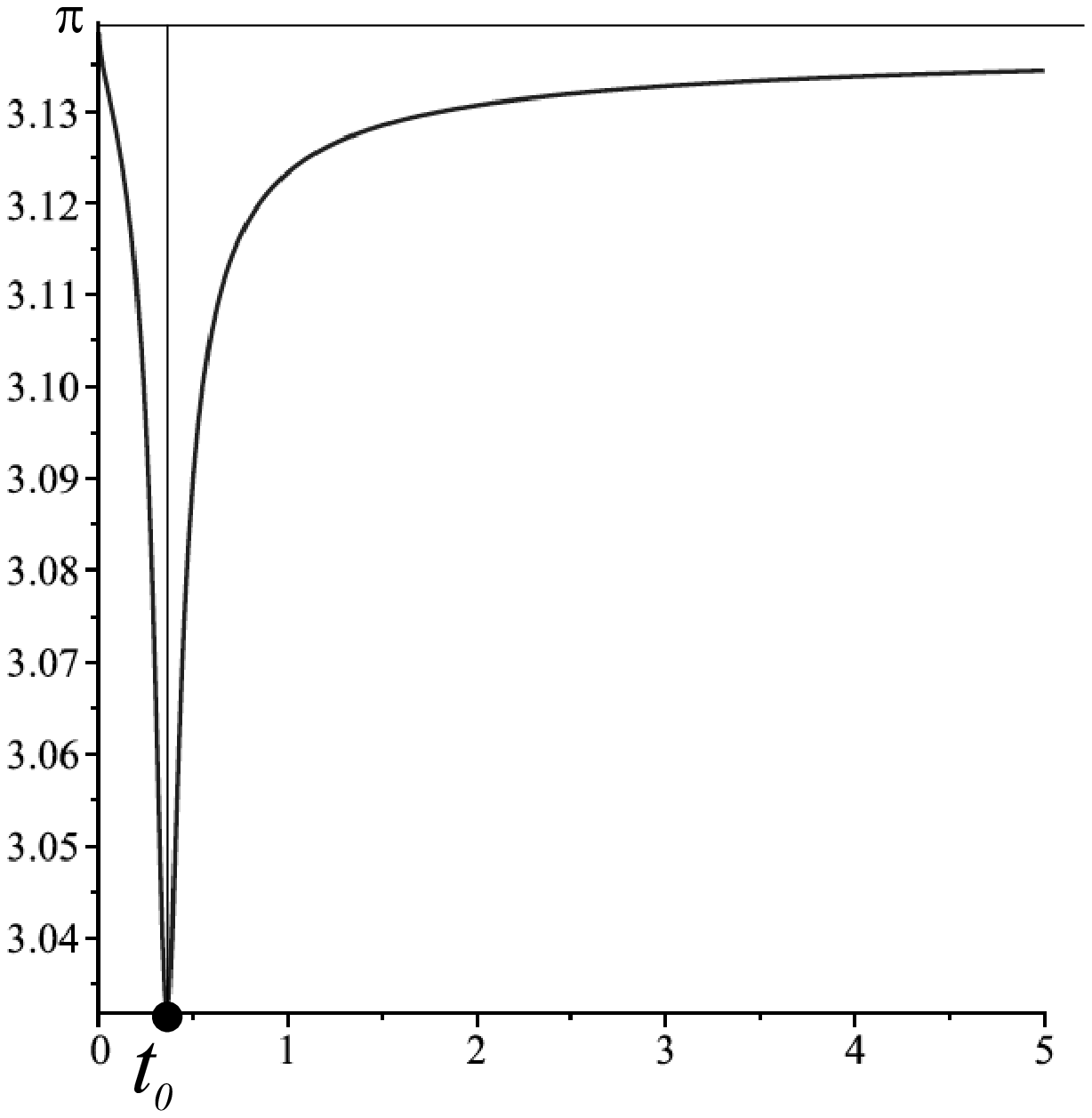}
\caption{$S(\Delta^{\HXR}(t))$ function related to parameters $x_2=2, y_2=3/1, z_2=1$ $x_3=3\cdot t,y_3=-1 \cdot t,z_3=0$.}
\label{}
\end{figure}
In Fig.~10 we described the $S(\Delta(t))$ function related to geodesic triangle $\Delta(t)$ $(t \in (0,5))$ with vertices $A_1=(1,1,0,0)$, $A_2=(1,2,3/2,1)$, $A_3=(1,3\cdot t,-1 \cdot t,0)$.
Its minimum is achieved at $t_0 \approx 0.36392$ where $S(\Delta^{\HXR}(t_0))\approx 3.03236$. ~ ~ ~ $\square$

Finally we obtain the following
\begin{thm}
The sum of the interior angles of a geodesic triangle of $\HXR$ space is less or equal to $\pi$. ~ ~ $\square$
\end{thm}
\medbreak

\end{document}